# TAILOR-MADE TESTS FOR GOODNESS OF FIT TO SEMIPARAMETRIC HYPOTHESES

By Peter J. Bickel,[1] Ya'acov Ritov[2] and Thomas M. Stoker

*University of California at Berkeley, Hebrew University of Jerusalem and Massachusetts Institute of Technology*

We introduce a new framework for constructing tests of general semiparametric hypotheses which have nontrivial power on the $n^{-1/2}$ scale in every direction, and can be tailored to put substantial power on alternatives of importance. The approach is based on combining test statistics based on stochastic processes of score statistics with bootstrap critical values.

**1. Introduction.** The practice of statistical testing plays several roles in empirical research. These roles range from the careful assessment of the evidence against specific scientific hypotheses to the judgment of whether an estimated model displays decent goodness of fit to the empirical data. The paradigmatic situation we consider is one where the investigator views some departures from the hypothesized model as being of primary importance, with others of interest if sufficiently gross, but otherwise secondary. For instance, low-frequency departures from a signal hypothesized to be constant might be considered of interest, even if of low amplitude; while high-frequency departures are less so, unless they are of high amplitude.

The optimal testing of a simple hypothesis against a simple alternative is the cornerstone of modern statistical theory. However, there is no clear notion of optimality for more complicated situations. The Hájek–Le Cam asymptotic theory proved that there exist strong concepts of asymptotic efficiency in parametric estimation. These ideas have been extended to semiparametric models—see [3, 14, 22]. However, there is no compelling sense of an asymptotically optimum test, in either the parametric or the semiparametric asymptotic theories, save for some simple one-parameter hypotheses.

Received November 2001; revised June 2005.
[1]Supported in part by NSF Grant FD-01-04075.
[2]Supported in part by ISF Grant 793.03.

We deal exclusively with the "elementary" case of i.i.d. data for ease of exposition. Moreover, all our considerations are asymptotic save for illustrative simulations. Generalization of this point of view to the two-sample problem, independent nonidentically distributed case, time series, and so on, is conceptually not difficult, and may be even simpler because of availability of permutation tests.

The general types of tests that people have constructed fall into one of two classes:

(i) Those which have nonnegligible asymptotic power for departures on the $n^{-1/2}$ scale in every possible direction. In the standard problems of testing goodness of fit to a single distribution against all alternatives, these are the classical tests of Kolmogorov and Cramér–von Mises and their classical extensions to the problem of testing fit to a parametric hypothesis on the one hand and independence on the other.

(ii) Those which have trivial asymptotic power at the $n^{-1/2}$ scale in every direction. The $\chi^2$ tests with increasing number of cells as $n \to \infty$ are the preeminent example of this type, but a number of variants have recently been explored through devices such as empirical likelihood—see [10] for recent examples.

Tests of type (ii) have the feature that they have approximately equal power in all directions. As a consequence, they can enjoy minimax properties over suitable nonparametric families of alternatives—see [15], for example. But, as we noted, they pay for this by not having power at rate $n^{-1/2}$ in any particular direction. The tests of type (i) have the weakness that they concentrate their power at the $n^{-1/2}$ scale in very explicit alternative directions, dictated primarily by the metric, implicitly or explicitly, used. For example, the Kolmogorov test for goodness of fit to the uniform $(0, 1)$ distribution is well known to have power mainly against alternatives such that $|P(X \leq 1/2) - 1/2|$ is large.

The principal reason for limiting oneself to tests of types (i) and (ii) appears to have been the need for simple approximations to the critical values under the null which need to be coupled with specification of a test statistic to implement a test. However, the critical values can be approximated by bootstrap methods, as discussed in Section 3.3.

Our goal in this paper is to show that it is possible to construct tests for any semiparametric hypothesis, which have as much power as possible at the $n^{-1/2}$ scale in a few directions of interest, specific to the particular scientific problem investigated, reserving some power for gross departures (in the $n^{-1/2}$ scale) for other directions.

We clearly do not adopt the minimax and adaptive minimax testing point of view of Ingster [15]. Our proposal does not aim at minimaxity and since we concentrate on the $n^{-1/2}$ scale our tests do not have uniformity properties except over relatively small families. We believe that in testing, even more than in estimation, prior information or biases need to be paid attention to, since, as Janssen [16] points out, achieving reasonable power over more than a few orthogonal directions is hopeless.

There has been another direction that we want to mention but do not develop in this paper. The idea is to construct a sequence of tests which are consistent against broader and broader classes of alternatives as one proceeds down the sequence, stopping testing at a data-determined point on the sequence. A limited proposal of this type was made by Rayner and Best [24] and developed more generally in [6]. Some important special cases are discussed in [13], Chapter 7, in the context of testing the hypothesis of no effect in nonparametric regression.

Our general approach, which is detailed in Section 3, is to use as building blocks one-dimensional score (Rao) test statistics for simple hypotheses. For composite hypotheses we use the natural generalizations of Rao tests, efficient scores. These efficient score tests are called Neyman $C(\alpha)$ tests in the statistics literature or conditional moment restriction tests in the econometrics literature (see [2]).

Conceptually, as we discuss in Section 3, our approach applies to general semiparametric hypotheses such as independence, the Cox model in survival analysis and index models in econometrics. It also, as we demonstrate, guides us how to proceed when we test a parametric or semiparametric model within a semiparametric alternative, for instance, independence within copula models, simple index versus multiple index models. We view this as the most important nontechnical contribution of the paper.

In Section 3 we give some general conditions under which the asymptotic theory for the types of test statistics discussed in this section, and for appropriate bootstrap critical values, is valid. We also study the power behavior of these tests under these assumptions. In Section 4 we discuss the classical examples of goodness of fit to a parametric hypothesis and independence. We show how the classical tests of Kolmogorov–Smirnov and Cramér–von Mises type fit into our framework, and also derive a variety of new tests based on our principles. We indicate how the general conditions of Section 3 are implied by mild and easily checkable conditions in these classical situations. We have chosen to exhibit the approach in detail in two situations here, namely parametric hypotheses and independence. Tests of index models are covered in [7]. However, as we have indicated, our approach is applicable to any of the hundreds of goodness-of-fit problems that arise with semiparametric models.

## 2. Heuristics.

2.1. *Parametric tests.* Suppose that $X_1, \ldots, X_n$ are a (i.i.d.) random sample from the probability $P \in \mathcal{Q}$, where $P \ll \mu$, $p = \frac{dP}{d\mu}$. Suppose that $\mathcal{Q} = \{P_\theta : \theta \in \mathbf{R}\}$ is a regular (one-dimensional) submodel of probabilities. Consider testing the hypothesis $H : P = P_0$ against $K : P = P_\theta$ where $\theta > 0$.

Denote the log-likelihood of an observation and its derivative at $\theta = 0$, the efficient score (see [23]) by

$$
(1) \qquad \ell(\cdot, P) \equiv \ln p, \qquad \dot{\ell}(\cdot, P_0) \equiv \frac{\partial \ell(\cdot, P_0)}{\partial \theta},
$$

where $\dot{\ell} \in L_2(P_0)$, $E_0(\dot{\ell}(\cdot, P_0)) = 0$. The familiar scoring test of Rao [23] is based on the mean score test statistic, and is locally asymptotically most powerful. Since $n^{-1/2} \sum_{i=1}^n \dot{\ell}(X_i, P_0) \xrightarrow{D} \mathcal{N}(0, \|\dot{\ell}\|_0^2)$ under the null, it uses as asymptotic critical value $z_{1-\alpha} \|\dot{\ell}\|_0$, where $\|h\|_0^2 = \int h^2 \, dP_0$ and $z_{1-\alpha}$ is the standard Gaussian $1 - \alpha$ quantile. Note that this test is consistent if and only if $E_\theta(\dot{\ell}(X, P_0)) > 0$ for $\theta > 0$, namely if a nonzero $\theta$ implies a positive mean score.

All the nontechnical difficulties in testing are already in composite parametric hypotheses, although they are traditionally ignored. Suppose that the general family is defined as $\mathcal{P} = \{P_{(\eta,\theta)} : \theta \in \mathbf{R}^q, \eta \in \mathbf{R}^p\}$ so that the general log-density takes the form $\ell = \ell(\cdot, P_{(\eta,\theta)})$, which for simplicity we write $\ell(\eta, \theta)$.

Suppose first that $q = 1 < p$. The null hypothesis is the restricted family $\mathcal{P}_0 \equiv \{P_{(\eta,0)}, \eta \in \mathbf{R}^p\}$. The null set can be approximated locally by the tangent space associated with the null hypothesis,

$$
(2) \qquad \overset{\bullet}{\mathcal{P}}_0(P_0) = \operatorname{span}\left\{ \frac{\partial \ell(\eta_0, 0)}{\partial \eta_j} : 1 \leq j \leq p \right\}.
$$

The score function of the alternative is $\dot{\ell} \equiv \partial \ell(\eta_0, 0)/\partial \theta$. However, part of this score is in fact in the null space. Therefore, the efficient score is that part that is left by removing contributions from directions in the tangent space,

$$
\ell^* = \frac{\partial \ell(\eta_0, 0)}{\partial \theta} - \sum_{j=1}^p a_j(\eta, 0_0) \frac{\partial \ell(\eta_0, 0)}{\partial \eta_j},
$$

where the $a_j$'s are projection (least squares) weights; namely $\{a_j(\eta_0, 0)\}$ minimize $\|\dot{\ell} - \sum_{j=1}^p a_j(\eta_0, 0) \, \partial \ell(\eta_0, 0)/\partial \eta_j\|_0^2$. If $q = 1$ and $\hat{\eta}$ is a $\sqrt{n}$-consistent estimator of $\eta$ under $H$, then the Neyman [21] $C(\alpha)$ test statistic is $T = n^{-1/2} \sum_{i=1}^n \ell^*(X_i, \hat{\eta}, 0)$.

When $q > 1 = p$, $\mathcal{P} = \{P_\theta : \theta \in \mathbf{R}^q\}$, the departure can happen in different directions. The tangent space $\overset{\bullet}{\mathcal{P}}(P_0)$ is the linear closure of $\{\partial \ell/\partial \theta_j : j =$

$1, \ldots, q\}$, and a standard test statistic is

$$T = \frac{1}{n}\left[\sum_{i=1}^{n}\nabla_\vartheta \ell(X_i)\right]' \mathcal{I}_0^{-1}\left[\sum_{i=1}^{n}\nabla_\vartheta \ell(X_i)\right],$$

where $\mathcal{I}_0$ is the information matrix. The Rao tests, which are called *Lagrange multiplier* tests in econometrics, have the advantage of making use of estimates of the statistical model only under the null hypothesis. In contrast, Wald tests [27] and likelihood ratio tests are based on comparing estimates of the model under alternatives with those of the model estimated under the null. In the parametric context, to first order, both Wald and likelihood ratio tests are equivalent to score-based tests.

It is important to note that $T$ is just one way of combining the different test directions. *There is nothing magic in the Mahalanobis distance.* Suppose that we can rank the alternative one-dimensional models for which $h_1, \ldots, h_q$ are score functions in order of plausibility. If they are orthogonal, it is plausible to use $T_a \equiv \sum_{j=1}^{q} \lambda_j^2 (n^{-1/2} \sum_{i=1}^{n} h_j(X_i))^2$ where $0 < \lambda_1 < \cdots < \lambda_p$ reflect the relative importance of the $h_j$. In general we would arrive at

$$(3) \qquad T_a = \frac{1}{n}\left[\sum_{i=1}^{n}\tilde{h}(X_i)\right]' \Lambda \Sigma_0^{-1} \Lambda \left[\sum_{i=1}^{n}\tilde{h}(X_i)\right],$$

where $\Lambda = \operatorname{Diag}(\lambda_1, \ldots, \lambda_p)$.

Another alternative is to use the union-intersection principle of Roy [25], to obtain

$$(4) \qquad T_b = \max_{1 \le j \le p}\left|\frac{\lambda_j}{\sqrt{n}}\sum_{i=1}^{n}h_j(X_i)\right|.$$

In general, any norm of the vector $(h_1, \ldots, h_q)$ could be used as a test statistic.

2.2. *Semiparametric essentials.* When the hypotheses are composite and semiparametric, the collection of score functions that generalize (1) and (2) depends on $P$, and no longer consists of Euclidean spaces.

Let $\mathcal{Q} = \{P_\theta : \theta \in \mathbf{R}\} \subset \mathcal{P}$ denote a regular one-parameter submodel with $P_0$ the true distribution. Clearly, the score in the model $\mathcal{Q}$, $h_\mathcal{Q} \equiv \dot{\ell}$ depends on $\mathcal{Q}$. A composite hypothesis $\mathcal{P}$ is the union of many, usually an infinite number of such $\mathcal{Q}$'s. The relevant set of (alternative) scores is the tangent space $\overset{\bullet}{\mathcal{P}}(P_0)$ defined as the linear closure of all the associated scores $h_\mathcal{Q}$—see, for instance, [3], Chapter 2, for details.

We parameterize the model by $\mathcal{P} = \{P_{(\alpha,\beta)} : \alpha \in A, \beta \in B\}$ where $A$, $B$ are subsets of function spaces, and the null hypothesis is $H : \beta = 0$. Define the

"full," "null" and "alternative" tangent spaces (see [3], page 70):

$$\overset{\bullet}{\mathcal{P}}(\alpha,\beta) = \text{tangent space of the model } \mathcal{P} \text{ at } P_{(\alpha,\beta)},$$

$$\overset{\bullet}{\mathcal{P}}(\alpha,0) = \text{tangent space of } \{P_{(\alpha,0)} : \alpha \in A\} \text{ at } P_{(\alpha,0)},$$

$$\overset{\bullet}{\mathcal{P}}_0^\perp(\alpha,0) = \text{orthogonal complement of } \overset{\bullet}{\mathcal{P}}_0(\alpha,0) \text{ in } \overset{\bullet}{\mathcal{P}}(\alpha,0).$$

That is, $\overset{\bullet}{\mathcal{P}}_0(\alpha,0) \perp \overset{\bullet}{\mathcal{P}}_0^\perp(\alpha,0)$ and $\overset{\bullet}{\mathcal{P}}(\alpha,0) = \overset{\bullet}{\mathcal{P}}(\alpha,0) \oplus \overset{\bullet}{\mathcal{P}}_0^\perp(\alpha,0)$.

The space $\overset{\bullet}{\mathcal{P}}_0(P_0)$ captures the directions of variation from $P_0$ that are consistent with the null hypothesis of interest. To test $\theta = 0$, we should remove from any $\dot{\ell} \in \overset{\bullet}{\mathcal{P}}(P_0)$ its component that is actually consistent with the null hypothesis and is in $\mathcal{P}_0(P_0)$. Therefore, the effective direction of interest for the alternative $\mathcal{Q}$ is given by *efficient score function*

$$\ell^*(\cdot,P_0) \equiv \dot{\ell}(\cdot,P_0) - \Pi(\dot{\ell},\alpha), \qquad \dot{\ell}(\cdot,P_0) \in \overset{\bullet}{\mathcal{P}}(P_0),$$

where $\Pi(\cdot,\alpha)$ is the projection operator from $L_2(P_\alpha)$ to the subspace $\overset{\bullet}{\mathcal{P}}_0(P_\alpha)$ of $L_2^0(P_\alpha)$, the space of square integrable functions with mean zero under $P_\alpha$.

NOTATION. We write $Q(h)$ for $\int h\, dQ$, and $P_n$ for the empirical distribution function.

### 3. Score tests.

3.1. *The score process and the testing paradigm.* The above ideas motivate our general testing paradigm.

Let $\Gamma$ be some index set, and let $h \equiv h(\gamma,\alpha) \in L_2^0(P_\alpha)$, $\gamma \in \Gamma$, $\alpha \in A$, be some test function. Recall that $A$ is the parameter space under the null. The *score process* is

$$Z_n(\gamma,\alpha) \equiv \frac{1}{\sqrt{n}} \sum_{i=1}^n \Pi^\perp(h,\alpha)(X_i),$$

where $\Pi^\perp$ is the projection from $L_2^0(P_\alpha)$ to $\overset{\bullet}{\mathcal{P}}_0^\perp(\alpha,0)$.

$\Gamma$ is an index set, pointing to a direction $h(\gamma,\alpha)$ in the tangent space, or more generally in $L_2^0(P_0)$, such that $\{h(\cdot,\alpha)\}$ is not too big, say a universal Donsker class. As we shall see in examples, the reason for making $h$ depend on $\alpha$ also is that it is natural to have the family of scores depend on where we think we are in $A$. To avoid technicalities, we assume that for all $x$ and $\alpha$, $h_\gamma(x,\alpha) \in l_\infty(\Gamma)$, the space of all bounded real-valued functions on $\Gamma$. We may write $h_\gamma(\alpha)$, suppressing dependence on $x$.

In general, $Z_n(\gamma, \alpha)$ is not computable given the data, but if $\hat{\alpha} \in A$ is an estimate of $\alpha$ we can consider

$$\hat{Z}_n(\gamma) \equiv Z_n(\gamma, \hat{\alpha})$$

defined on $\Gamma$. If $\hat{\alpha}$ is an MLE, $\hat{Z}_n$ simplifies as in the parametric case to

$$\hat{\hat{Z}}_n(\gamma) = \sqrt{n}(P_n - P_{\hat{\alpha}})(h_\gamma(\cdot, \hat{\alpha}))$$

since $P_n(v) = 0$ for all $v \in \overset{\bullet}{\mathcal{P}}_0(\hat{\alpha})$ is a restatement of the likelihood equations. In particular, $P_n(\Pi(h, \hat{\alpha})) = 0$. We will also consider $\hat{\hat{Z}}_n$ more generally for $\hat{\alpha}$ an efficient estimate in the sense of [3], Chapter 5, pages 179–182.

We think of $\hat{Z}_n(\cdot)$, $\hat{\hat{Z}}_n(\cdot)$, and so on as stochastic processes defined on $\Gamma$ related to empirical processes—see [26], for instance. We shall use $\hat{Z}_n$ and $\hat{\hat{Z}}_n$ to construct tailor-made tests.

Let $A_0 \subset A$ be a neighborhood of the true $\alpha_0$ where $A$ is a metric space with metric $\rho$. As above, we write $P_0$ for $P_{\alpha_0}$, and so on. We always require

(5) $$Z_n(\cdot, \alpha_0) \Rightarrow Z(\cdot, \alpha_0)$$

under $P_{\alpha_0}$, for all $\alpha_0 \in A_0$, in the sense of weak convergence for $l^\infty(\Gamma)$-valued variables, where $Z(\cdot, \alpha_0)$ is a mean-zero Gaussian process with

$$\text{cov}(Z(\gamma_1, \alpha_0), Z(\gamma_2, \alpha_0)) = \text{cov}_0(h_1, h_2) - \text{cov}_0(\Pi(h_1, \alpha_0), \Pi(h_2, \alpha_0))$$

with the obvious convention in notation. (To be exact, we should be speaking of outer probabilities since we interpret weak convergence in the sense of Hoffman–Jørgensen. But measurability issues can be dealt with easily in the situations we are interested in and we ignore them in the future.) The property we want is

(6) $$\hat{\hat{Z}}_n \quad \text{or} \quad \hat{Z}_n(\cdot) \Rightarrow Z(\cdot, \alpha_0)$$

in the same sense as above.

We propose to base our tests on the score process, at least conceptually. What (6) will give us is the weak convergence of statistics of the form $T(\hat{Z}_n(\cdot))$ where $T : l_\infty(\Gamma) \to R$ continuously. Possibilities are

(7) $$T_\mu \equiv \int \hat{Z}_n^2(\gamma) \, d\mu(\gamma)$$

for $\mu$ a finite measure on $T$, or

(8) $$T_K \equiv \sup_\Gamma |\hat{Z}_n(\gamma)|$$

or more general $\mu$ norms of $|\hat{Z}_n|$, or even $\alpha$-dependent $\mu$'s which are suitably continuous in $\alpha$. By taking the span $\{h(\gamma, \alpha) : \gamma \in \Gamma\}$ dense in $\overset{\bullet}{\mathcal{P}}(\alpha, 0)$ for all $\alpha$, and $\mu$ with support $\Gamma$, we can expect consistency against all alternatives. We will illustrate further in the examples of the next section. A simple example is the following.

EXAMPLE 3.1 (Goodness-of-fit statistics). Consider testing the null hypothesis that a distribution on $\mathbf{R}$ is $P_0$ against "all" alternatives, namely where $\overset{\bullet}{\mathcal{P}}(P_0) = L_2^0(P_0)$. We consider the family of directions $h_\gamma(\cdot) = 1(\cdot \leq \gamma) - F_0(\gamma)$, $\gamma \in \mathbf{R}$, where $F_0$ is the cumulative distribution function of $P_0$. The following two statistics arise in association with those above. Associated with (3) is the familiar *Cramér–von Mises* (CvM) goodness-of-fit statistic, $T_a = n \int (F_n(\gamma) - F_0(\gamma))^2 \, dF_0(\gamma)$, where $F_n$ is the empirical distribution function, and the weighting measure is $\mu = F_0$. Corresponding to (4) is the familiar *Kolmogorov–Smirnov* (KS) goodness-of-fit test statistic, $T_b = \sup_\gamma |\sqrt{n}(F_n(\gamma) - F_0(\gamma))|$. Note that these $h_\gamma$ here are typically chosen because we start with the cumulative distribution function as our representation of the probability $P$, not because of a desire for power in clearly defined directions (which is the result).

3.2. *General theorems.* We close this section with some general theorems giving essentially the minimal conditions under which our heuristics for test statistic construction and critical value setting are justified. Checking the conditions of these theorems is the major difficulty.

Here are the conditions we use for our theorems. The estimate $\hat{\alpha}$ is such that for $P \equiv P_0 \equiv P_{\alpha_0}$:

(M0) $\{h_\gamma(\cdot, \alpha_0) - \Pi(h_\gamma(\cdot, \alpha_0), \alpha_0) : \gamma \in \Gamma\} = \{\Pi^\perp(h_\gamma(\cdot, \alpha_0)) : \gamma \in \Gamma\}$ is a universal Donsker class.
(M1) $\|(P_n - P_0)(\Pi(h_\gamma(\cdot, \hat{\alpha}), \hat{\alpha}) - \Pi(h_\gamma(\cdot, \hat{\alpha}), \alpha_0))\|_\infty = o_P(n^{-1/2})$.
(M2) $\sup\{\|(P_{\hat{\alpha}} - P_0)h_\gamma(\cdot, \alpha) + P_0 \Pi(h_\gamma(\cdot, \alpha), \hat{\alpha})\|_\infty : \alpha \in A\} = o_P(n^{-1/2})$.
(M3) $\|(P_n - P_0)(h_\gamma(\cdot, \hat{\alpha}) - h_\gamma(\cdot, \alpha_0))\|_\infty = o_P(n^{-1/2})$.
(M4) $\|(P_{\hat{\alpha}} - P_0)(h_\gamma(\cdot, \hat{\alpha}) - h_\gamma(\cdot, \alpha_0))\|_\infty = o_P(n^{-1/2})$.
(M5) $\|(P_{\hat{\alpha}} - P_0)h_\gamma(\cdot, \alpha_0) - (P_n - P_0)\Pi(h_\gamma(\cdot, \alpha_0), \alpha_0)\|_\infty = o_P(n^{-1/2})$.

NOTES. 1. (M3) and (M4) are automatically satisfied if $h_\gamma(\cdot, \alpha)$ does not depend on $\alpha$.

2. If $\mathcal{H}_0 \equiv \{h_\gamma(\cdot, \alpha_0) - P_{\alpha_0} h_\gamma(\cdot, \alpha_0) : \gamma \in \Gamma\}$ is a Donsker class, showing that $\{\Pi(h, \alpha_0) : h \in \mathcal{H}_0\}$ is also Donsker is usually not hard. For instance, suppose $\Pi$ preserves order, $\Pi(h_1, \alpha_0) \leq \Pi(h_2, \alpha_0)$ if $h_1 \leq h_2$. If $\mathcal{H} \equiv \{h_\gamma(\cdot, \alpha_0) : \gamma \in \Gamma\}$ satisfies the bracketing entropy condition of Theorem 2.8.4, page 172, of [26], then since $\Pi(\cdot, \alpha_0)$ is $L_2(P_0)$ norm reducing, $\{\Pi(h, \alpha_0) : h \in \mathcal{H}_0\}$ also satisfies the same condition. Thus (M0) is usually not difficult.

3. (M5) says that $\hat{\alpha}$ is efficient under $H$, a generalization of the requirement that $\hat{\alpha}$ be (a regularly behaving) MLE—see [3], pages 176–182.

Here are two theorems.

THEOREM 3.1. *If* (M0), (M3), (M4) *and* (M5) *hold, then for all* $\alpha_0$,

(9) $$Z_n(\cdot, \alpha_0) \Rightarrow Z(\cdot, \alpha_0),$$

(10) $$\hat{\hat{Z}}_n(\cdot) = Z_n(\cdot, \alpha_0) + o_{P_0}(1),$$

*and hence,*

(11) $$\hat{\hat{Z}}_n(\cdot) \Rightarrow Z(\cdot, \alpha_0).$$

PROOF. By construction and (M3)

$$\hat{\hat{Z}}_n(\gamma) = n^{1/2}\{(P_n - P_0)(h_\gamma(\cdot, \hat{\alpha})) - (P_{\hat{\alpha}} - P_0)(h_\gamma(\cdot, \hat{\alpha}))\}$$
$$= n^{1/2}\{(P_n - P_0)(h_\gamma(\cdot, \alpha_0)) - (P_{\hat{\alpha}} - P_0)(h_\gamma(\cdot, \hat{\alpha}))\} + o_P(1)$$
$$= n^{1/2}\{(P_n - P_0)(h_\gamma(\cdot, \alpha_0)) - (P_{\hat{\alpha}} - P_0)(h_\gamma(\cdot, \alpha_0))\} + o_P(1)$$

by (M4). Finally by (M5),

$$\hat{\hat{Z}}_n(\gamma) = n^{1/2}\{(P_n - P_0)(h_\gamma(\cdot, \alpha_0) - \Pi(h_\gamma(\cdot, \alpha_0), \alpha_0))\} + o_P(1),$$

which is just (10). Note that $o_P(1)$ is interpreted here in the sense of $\|\cdot\|_\infty$ on functions of $\gamma$. Conclusions (9) and (11) follow immediately from (M0) and (10). □

THEOREM 3.2. *If* (M0)–(M3) *hold, then*

$$\hat{Z}_n(\gamma) = Z_n(\gamma) + o_{P_0}(1).$$

The proof appears in the Appendix.

Condition (M2) is implied by the following two more easily checkable conditions. See Lemma A.1.

(N1) (i) $\hat{\alpha}$ is consistent in the Hellinger metric $\rho_H$ given by $\rho_H^2(\alpha_1, \alpha_2) = \int(\sqrt{dP_{\alpha_1}} - \sqrt{dP_{\alpha_2}})^2$.

(ii) Let $A_0$ be a fixed Hellinger ball around $\alpha_0$ and suppose that $\mu \gg P_\alpha$, $\alpha \in A_0$, and with $\mu$ a probability measure. Let $\|\cdot\|_\mu$ be the $L_2(\mu)$ norm. Write $s(\alpha) = \sqrt{dP_\alpha/d\mu}$, $\hat{s} \equiv s(\hat{\alpha})$, $s_0 \equiv s(\alpha_0)$, assume

$$\|(\hat{s} - s_0)^2/\hat{s}\|_\mu^2 = \int (s_0/\hat{s} - 1)^4 \hat{s}^2 \, d\mu = o_P(n^{-1}).$$

(iii) Let $\Pi_\mu$ denote projection in $L_2(\mu)$ onto the tangent space at $s(\alpha_0)$ of $\mathcal{L} = \{s(\alpha): \alpha \in A\}$. Assume, $\|\hat{s} - s_0 - \Pi_\mu(\hat{s} - s_0)\|_\mu = o_P(n^{-1/2})$.

(N2) (i) $\sup\{\|h_\gamma(\cdot,\alpha)\|_\infty : \gamma \in \Gamma, \alpha \in A_0\} < \infty$.
(ii) $\sup\{\|\Pi(h_\gamma(\cdot,\alpha))\|_\infty : \gamma \in \Gamma, \alpha \in A_0\} < \infty$.

NOTE. In smooth parametric models (N1)(iii) holds if $\|\hat{s} - s_0\|_\mu = o_P(n^{-1/4})$, $\|s - s_0 - \Pi_\mu(s - s_0)\|_\mu = O(\|s - s_0\|_\mu^2)$.

We may wish to consider (see below) statistics in which the averaging measure also depends on $\alpha$, say $\int \hat{Z}_n^2(h_\gamma) \, d\mu(\gamma, \hat{\alpha})$. This too can be dealt with by a condition such as $\alpha \to \mu(\cdot, \alpha)$ is uniformly continuous on $\mathcal{H}$ in the bounded variation topology on the finite signed measures on $\Gamma$. More generally, we may simply consider any test statistic of the form $F(\hat{Z}_n, \hat{\alpha})$, where $F: l_\infty(\Gamma) \times A \to R$ is continuous in the $l_\infty \times \rho$ topology.

3.3. *Setting critical values*: *the bootstrap.* There is a novel issue that arises in the context of composite hypotheses. The statistic $\hat{Z}_n \equiv \frac{1}{\sqrt{n}} \sum_{i=1}^n l^*(X_i, P_{\hat{\alpha}})$ arising from the situation where there is only one direction of departure is, if $P_\alpha$ is true, an approximation to $Z_n(\alpha) \equiv \frac{1}{\sqrt{n}} \sum_{i=1}^n l^*(X_i, P_\alpha)$. Since $Z_n(\alpha)$ has an $\mathcal{N}(0, I(P_\alpha))$ limiting distribution, a Gaussian critical value using $I(P_{\hat{\alpha}})$ is appropriate. On the other hand, if the null hypothesis is simple, critical values for any statistic can be obtained by simulation. But in the general situation of composite hypotheses, that we now consider, unless there is invariance, the most plausible way of setting critical values is by a bootstrap. The natural choice is to simulate from $P_{\hat{\alpha}}$. That is, let $S(\gamma, \alpha) = \Pi^\perp(h, \alpha)$, where $h \equiv h(\gamma, \alpha)$, and let

$$\tilde{\hat{Z}}_n(\gamma) = n^{-1/2} \sum_{i=1}^n S(\gamma, \hat{\alpha})(\tilde{X}_i),$$

where the $\tilde{X}_i$'s are i.i.d. from $P_{\hat{\alpha}}$. We expect that if (6) holds, $\tilde{\hat{Z}}_n \Rightarrow Z(\cdot, \alpha_0)$ in $P_0$ probability. That is, the Prohorov distance between the $P_{\hat{\alpha}}$ distribution of $\tilde{\hat{Z}}_n$ and the distribution of $\hat{Z}_\alpha$ tends to 0 in $P_\alpha$ probability.

To ensure that this bootstrap method works for both $\hat{Z}_n$ and $\tilde{\hat{Z}}_n$ we need to simply replace conditions (M0)–(M5) by versions uniform in $\alpha_0 \in A_0$. We leave a formal statement to the reader. Alternatively, in these cases, as has been explored in [2] and [4], it is also possible to use the $m$ out of $n$ bootstrap, simulating the distribution of the statistic for samples of size $m$ by drawing subsamples of size $m$ from the original sample, where $m \to \infty$, $m/n \to 0$.

There is another way of bootstrapping discussed in [5] which may be simpler since it only involves resampling. Let

$$\hat{Z}_n^*(\cdot) = n^{-1/2} \sum_{i=1}^n (S(\gamma, \hat{\alpha}^*)(X_i^*) - S(\gamma, \hat{\alpha})(X_i)),$$

where $\hat{\alpha}^* \equiv \hat{\alpha}(X_1^*, \ldots, X_n^*)$ and $X_1^*, \ldots, X_n^*$ are i.i.d. from the empirical distribution $P_n \equiv n^{-1} \sum_{i=1}^n \delta_{X_i}$. The appropriate heuristic is that again if (6) holds, $\hat{Z}_n^* \Rightarrow Z(\cdot, \alpha_0)$ in $P_0$ probability.

This bootstrapping method is more problematic to check. Essentially what is needed for $\hat{Z}_n^*$ to obey (3.5) are conditions given in [5]. For these we make the following identification: Suppose $\hat{\alpha} = \alpha(P_n)$ where $\alpha : \mathcal{M} \to A$, $\alpha(P_\alpha) = \alpha$ for $P_\alpha \in \mathcal{P}$ and $\mathcal{M}$ is the set of all probabilities. Let $T : \mathcal{M} \to L(\mathcal{H}) \times L(\mathcal{H})$ be defined for $\mathcal{H}$ a Banach space of functions containing $\{h_\gamma : \gamma \in \Gamma\}$ and $L(\mathcal{H})$ the set of bounded linear functionals on $\mathcal{H}$ by

$$T(P)(h) = \left( \int \Pi(h, \alpha(P)) \, dP, \int h \, d(P - P_{\alpha(P)}) \right).$$

Note that $T(P_\alpha) \equiv 0$ so that the hypothesis is contained in $\{P : T(P) \equiv 0\}$. Now put on $T(P)$ the conditions specified by Bickel and Ren [5].

3.4. *Power.* It is easiest to see what happens to the processes on which we build our tests in the case where alternatives converge to $P_0$ in the $n^{-1/2}$ scale. Specifically suppose $\{P_t : |t| < 1\}$ is a one-dimensional regular parametric model through $P_0$ with score function $g(\cdot)$ such that $g \notin \overset{\bullet}{\mathcal{P}}_0(P_0)$, that is, it is possible to discriminate $\{P_t, t \neq 0\}$ from $\mathcal{P}_0$ at the $n^{-1/2}$ scale. Then, let $Z_g(h_\gamma, 0)$ be the Gaussian process with the same covariance structure as $Z(h_\gamma, 0)$ but with

$$E_0 Z_g(h_\gamma, 0) = \int (h_\gamma - \Pi(h_\gamma, \alpha_0)) g \, dP_0.$$

Evidently, $Z_g(h_\gamma, 0) = Z(h_\gamma, 0)$, if $g$ does not have a component orthogonal to $\overset{\bullet}{\mathcal{P}}_0$. Define $Z_g(h_\gamma, \alpha)$ similarly for $g \in L_2(P_\alpha)$, $g \notin \overset{\bullet}{\mathcal{P}}_0(P_\alpha)$. The following result is an immediate consequence of Le Cam's LAN theory; see, for example, his "third lemma" [12] and our theorems.

THEOREM 3.3. *Suppose $g \in \overset{\bullet}{\mathcal{P}}_0^\perp(P_\alpha)$ is, for each $\alpha$, the score function of a regular model through $P_\alpha$. Assume the sufficient conditions of Theorems* 3.1 *or* 3.2 *hold. Then $\hat{Z}_n(\cdot) \overset{t_n}{\Rightarrow} Z_g(\cdot, \alpha)$, where $\overset{t_n}{\Rightarrow}$ is weak convergence under $P_{g,t_n}$, where $\{P_{g,t} : |t| < 1\}$ is a regular model passing through $P_0 = P_\alpha$, with score function $g$ at $0$, and $t_n = t n^{-1/2}$ for fixed $t$.*

Suppose $q$ is bowl-shaped and symmetric and its discontinuity set has probability 0. That is, if $C \equiv \{z : q \text{ is continuous at } z\}$, $P[Z(\cdot, \alpha) \text{ or } Z_g(\cdot, \alpha) \notin C] = 0$, and $q : l_\infty(\Gamma) \to R$, $q(z) = q(-z)$, $q(\lambda z)$ strictly increasing in $\lambda$ for $\lambda > 0$ and all $z$. Then, if, as we assume, $Z(\cdot, \alpha)$ is tight, we have

(12) $$E(q(Z(\cdot, \alpha))) \leq E(q(Z_g(\cdot, \alpha))).$$

Equation (12) follows from Anderson's theorem if $\{G_1, \ldots, G_k\}$ forms a partition of $\Gamma$, $Z(\gamma, \alpha)$ is replaced by $Z^{(k)}(\gamma, \alpha) \equiv \sum_{j=1}^{k} Z(\gamma_j, \alpha) \mathbb{1}(\gamma \in G_j)$ and $Z_g$ is similarly approximated. Now $\|Z_g^{(k)}(\cdot, \alpha) - Z_g(\cdot, \alpha)\|_\infty \xrightarrow{P} 0$ as $k \to \infty$, for all $g$ including $g = 0$, and (12) follows in general.

It holds that test statistics of the form (7)–(8) have

$$\liminf_n P_{g,t_n}[T_n \geq c] > \liminf_n P_0[T_n \geq c]$$

for $t_n = \lambda n^{-1/2}$, all $\lambda > 0$, $c$ as desired, and for all $g$ such that $E_0 Z_g(h_\gamma, 0) \neq 0$ for some $\gamma$.

3.5. *Consistency.* Consistency against fixed $P \notin \mathcal{P}_0$ can be obtained by a strengthening of conditions—though the strengthening we now give is overkill.

Suppose that $P$ is as above. For a suitable $\alpha(P) \in A$ call $(M_P j)$, $j = 0, \ldots, 5$, condition $(Mj)$ with $P_0$ replaced by $P$ and $\alpha_0$ replaced by $\alpha(P)$. Define the process $Z_P(\cdot)$ as the Gaussian process with mean zero and the covariance structure given in (5) with $\Pi(\cdot, \alpha_0)$ replaced by $\Pi(\cdot, \alpha(P))$ and $P_0$ replaced by $P$. Then the conclusion of Theorems 3.1 and 3.2 holds if $P_0$ is replaced by $P$ with:

(i) $\hat{Z}_n(h_\gamma)$ replaced by $\hat{Z}_n(h_\gamma) - \sqrt{n} \int h_\gamma(\cdot, \alpha(P)) \, dP$;
(ii) $Z(\cdot, \alpha_0)$ replaced by $Z_P(\cdot)$.

We conclude that $|\hat{Z}_n(h_\gamma)| \xrightarrow{P} \infty$ if $\int h_\gamma(\cdot, \alpha(P)) \, dP \neq 0$. Thus consistency holds for $T_n$ given by (7) if $\int h_\gamma(\cdot, \alpha(P)) \, dP \neq 0$ for some $\gamma$ and all $P \neq P_0$. Consistency for other statistics can be reasoned analogously.

**4. Examples.** In this section we consider a few important examples in which we show how our notions produce tests which have appeared in the literature and some new ones. Our point is to illustrate the ideas of the score process and the tailor-made tests.

4.1. *Testing goodness of fit to a composite parametric hypothesis.* Let $\{P_{\boldsymbol{\theta}} : \boldsymbol{\theta} \in \Theta \subset R^d\}$ be a regular parametric model and let $\hat{\boldsymbol{\theta}}_n$ be a regular ([3], pages 18–19) estimate of $\boldsymbol{\theta}$ under $H$. We test $H$ against a saturated model such that $\overset{\bullet}{\mathcal{P}}(P) = \{a \in L_2(P) : Ea(X_1) = 0\}$:

$$\overset{\bullet}{\mathcal{P}}_0(P_\theta) = \operatorname{span}\left\{\frac{\partial l}{\partial \theta_1}(X_1, \boldsymbol{\theta}), \ldots, \frac{\partial l}{\partial \theta_d}(X_1, \boldsymbol{\theta})\right\},$$

where $l = \log p(x, \boldsymbol{\theta})$ is the log-likelihood. Then

$$\overset{\bullet}{\mathcal{P}}_0^\perp(P_\theta) = \left\{a(X_1) - \sum_{j=1}^d c_j(a, \boldsymbol{\theta}) \frac{\partial l}{\partial \theta_j}(X_1, \boldsymbol{\theta}) : a \in L_2(P_\theta), E_\theta a(X_1) = 0\right\}$$

and $c_j(a, \boldsymbol{\theta})$ is the coefficient of the projection of $a$ on $\overset{\bullet}{\mathcal{P}}_0(P_\theta)$, defined by minimizing

$$E_{\boldsymbol{\theta}}\left(a(X_1) - \sum_{j=1}^{d} c_j(a, \theta)\frac{\partial l}{\partial \theta_j}(X_1, \boldsymbol{\theta})\right)^2.$$

Identifying $\alpha$ with $\boldsymbol{\theta}$, $h = h_\gamma(\cdot, \boldsymbol{\theta})$, we obtain

$$S(\gamma, \boldsymbol{\theta}) = h - E_{\boldsymbol{\theta}} h(X_1) - \sum_{j=1}^{d} c_j(h, \boldsymbol{\theta}) \frac{\partial l}{\partial \theta_j}(X_1, \boldsymbol{\theta}).$$

The corresponding estimated score process is, for an estimate $\hat{\boldsymbol{\theta}}$, given by

$$\hat{Z}_n(h) = n^{-1/2} \sum_{i=1}^{n} \left\{ (h(X_i) - E_{\hat{\boldsymbol{\theta}}} h(X_i)) - \sum_{j=1}^{d} c_j(h; \hat{\boldsymbol{\theta}}) \frac{\partial l}{\partial \theta_j}(X_i, \hat{\boldsymbol{\theta}}) \right\}.$$

Suppose $\mathcal{P}_0$ is regular parametric, and more:

(R1) $\hat{\boldsymbol{\theta}}$ is regular on $\Theta$.
(R2) Suppose $\Gamma$ is compact $\subset \bar{R}^p$, where $\bar{R} = [-\infty, \infty]$, the processes $(\gamma, \theta) \to n^{1/2}(P_n - P_0)h_\gamma(\cdot, \theta)$ are tight and $\sup_{x, \gamma, \theta} |h_\gamma(x, \theta)| < \infty$.
(R3) The map $\theta \to h.(\cdot, \theta)$ is continuous in the norm on functions of $(\gamma, x)$ given by $\|\omega\|^2 = \sup_\gamma \int \omega^2(x, \gamma) \, dP_0(x)$.

The following proposition is a consequence of Theorem 3.2. We check its conditions using Lemma A.1. (M0) and (M1) hold and (N1)(i) and (N2) are immediate. We can check (M2) via (N1) and (N2). Condition (N1)(iii) follows since $\Pi_\mu(s(\hat{\theta}) - s(\theta_0)) = \overset{\bullet}{s}(\theta_0)(\hat{\theta} - \theta_0) + o_P(\hat{\theta} - \theta_0)$. Condition (N1)(ii) requires further conditions. For instance, it follows if the likelihood ratio $s(\cdot, \theta)/s(\cdot, \theta')$ is uniformly bounded for $\theta, \theta'$ within $\varepsilon$ of $\theta_0$, a case which unfortunately excludes the Gaussian, but the condition can be checked directly fairly easily for suitable $\hat{\theta}$. We have established

PROPOSITION 4.1. *If* (R1)–(R3) *and* (N1)(ii) *hold for a regular parametric hypothesis and given estimate* $\hat{\alpha}$, *then* $\hat{Z}_n(\cdot) \Rightarrow Z(\cdot)$ *where*

$$\operatorname{cov}(Z(h_1), Z(h_2)) = \operatorname{cov}_{\boldsymbol{\theta}_0}(h_1(X_1), h_2(X_1)) - \mathbf{c}^T(h_1, \boldsymbol{\theta}_0) I(\boldsymbol{\theta}_0) \mathbf{c}(h_2, \boldsymbol{\theta}_0),$$

*where* $\mathbf{c} = (c_1, \ldots, c_d)$ *and* $I$ *is the Fisher information.*

Versions of a result such as this one appear in [9] and [19] when the $h_\gamma$ are indicators of half lines.

The corresponding result for the bootstrap process $\tilde{\hat{Z}}_n$ is also valid as is that for $\hat{Z}_n^*$. That is, both the parametric and the Bickel–Ren application of the nonparametric bootstrap to testing can be used to set critical values.

Suppose $h$ does not depend on $\alpha$ and we weaken (R1) to

(R1)' $\hat{\theta} = \theta + o_{P_\theta}(n^{-1/4})$ for all $\theta \in \Theta$.

The theorem will still hold provided that we have Cramér conditions on several derivatives of the likelihood ensuring that the remainders in (N1) are quadratic in $\hat{\theta} - \theta_0$. Note that Proposition 4.1 enables us to plug in subefficient estimates without affecting the properties of our tests.

PROPOSITION 4.2. *If the hypothesis is regular parametric,* (R1), (R2) *and* (R3) *hold and $\hat{\theta}$ is efficient in the sense of* [3], *page* 43, *then the conclusion of Theorem* 3.1 *is valid.*

PROOF. We need only check (M5) and (M4). The former follows from

$$\sup_{\mathcal{H}} \left| \int h(x)p(x,\hat{\boldsymbol{\theta}})\,dx - \int h(x)p(x,\boldsymbol{\theta}_0)\,dx \right.$$
$$\left. - \int h(x)\dot{\ell}(x,\boldsymbol{\theta})(\hat{\boldsymbol{\theta}} - \boldsymbol{\theta}_0)p(x,\boldsymbol{\theta}_0)\,dx \right|$$
$$= o_P(|\hat{\boldsymbol{\theta}} - \boldsymbol{\theta}_0|),$$

and this requires in view of (R2) only that $\boldsymbol{\theta} \to p(\cdot, \boldsymbol{\theta})$ is $L_1$ differentiable, which is a consequence of regularity. The latter follows from regularity of $P_0$ and (R1) and (R3). □

Again with $\hat{\boldsymbol{\theta}}$, the MLE, results such as this one appear in [9] and [19]. The uniformity required for both versions of the bootstrap can easily be imposed.

4.2. *The Gaussian model.* We specialize to one of the most important parametric hypotheses $\mathcal{P} = \{\mathcal{N}(\mu, \sigma^2) : \mu \in R, \sigma^2 > 0\}$. Here we naturally take $\hat{\boldsymbol{\theta}} = (\bar{X}, \hat{\sigma}^2)$, the MLE's. It is convenient to use the invariance properties of the hypothesis and take $h(\cdot, \gamma, \boldsymbol{\theta}) = h_\gamma(\frac{x-\mu}{\sigma})$ if $\boldsymbol{\theta} = (\mu, \sigma)$. With this choice we are considering

$$\hat{\hat{Z}}_n(\gamma) = \frac{1}{\sqrt{n}} \sum_{i=1}^n \left( h_\gamma\left(\frac{X_i - \bar{X}}{\hat{\sigma}}\right) - \int_{-\infty}^\infty h(z)\varphi(z)\,dz \right).$$

If $\{h_\gamma = \mathbb{1}(-\infty, \gamma), \gamma \in R\} \equiv \mathcal{H}_0$, then satisfaction of (R1)–(R3) is easy. Using $\hat{\hat{Z}}_n$ as above we arrive at the common test statistics of Kac, Kiefer and Wolfowitz [17]: the *Kolmogorov–Smirnov* type,

$$\sup_{\mathcal{H}_0} |\hat{\hat{Z}}_n(\gamma)| = \sup_x \left| F_n\left(\frac{x-\bar{X}}{\hat{\sigma}}\right) - \Phi(x) \right|,$$

and the *Cramér–von Mises* type,

$$\int_{-\infty}^{\infty} \hat{\hat{Z}}_n^2(h_x)\,d\Phi(x) = \int_{-\infty}^{\infty} \left(F_n\left(\frac{x-\bar{X}}{\hat{\sigma}}\right) - \Phi(x)\right)^2 d\Phi(x),$$

where $F_n$ is the empirical d.f. $\mathcal{H}_0$ is a well-known universal Donsker class and the classical limiting result for Cramér–von Mises tests given in [17] follows. Tests can be implemented for both statistics using either of the two bootstraps. Invariance here implies that only simulation under $\mathcal{N}(0,1)$ is required. Other classes of tests are covered, for example, tests based on the empirical characteristic function [11].

We can also tailor statistics more carefully. For example, we can consider the two Gaussian mixture models as the alternative:

$$(1-\varepsilon)\Phi\left(\frac{t-\mu}{\sigma}\right) + \varepsilon\Phi\left(\frac{t-\mu-\Delta}{\sigma}\right), \qquad \mu,\Delta \in R, \sigma > 0, 0 < \varepsilon \leq \frac{1}{2}.$$

At least formally the tangent set $\overset{\bullet}{\mathcal{P}}(P_0)$ (see [3], page 50) at $\varepsilon = 0$, $\boldsymbol{\theta} = (\mu,\sigma)$ is just the set $\mathrm{span}\{X_1-\mu, (X_1-\mu)^2\} \cup \{(\exp\{\frac{\Delta}{\sigma^2}(X_1-\mu) - \frac{\Delta^2}{2\sigma^2}\} - 1): \Delta \in R\}$. We are led to consider $\mathcal{H}_0 = \{\exp\{\lambda x - \frac{\lambda^2}{2}\}: \lambda \in R\}$ and statistics such as

$$T_n \equiv \sup_\lambda \left|\frac{1}{\sqrt{n}}\sum_{i=1}^n (e^{\lambda((X_i-\bar{X})/\hat{\sigma})-\lambda^2/2} - 1)\right|.$$

Unfortunately, $\mathcal{H}_0$ is not a Donsker class and $T_n \overset{P}{\to} \infty$ under $H$; see [1]. Our heuristics and Theorem 3.1 apply if we restrict $\lambda$ to a compact set. The power against $n^{-1/2}$ alternatives of such $T_n$ persists. Note that $T_n$ can be viewed as a diagnostic since the maximizing value of $\lambda$ indicates where a second component might be.

We can also consider versions of the Cramér–von Mises approach reflecting our goals more precisely. For instance, consider a wavelet basis for $[0,1]$ written lexicographically $\omega_{ij}$ in order of scale and then within scale with $\omega_{11} \equiv 1$. Then given that we care more for departures at lower scales, consider $\lambda_{ij} = \lambda_i \sim \rho^i$, $1 \leq j \leq 2^i$, $\rho < \frac{1}{4}$. Since $\|\omega_{ij}\|_\infty = O(2^{i/2})$, if we let $h_{ij} = \omega_{ij}(\Phi(\cdot))$, and $h_t(x) \equiv \sum_{i,j} \lambda_{ij} h_{ij}(x) h_{ij}(t)$, then $\|h_t\|_\infty \leq M < \infty$ and $T = \int_0^1 \hat{\hat{Z}}_n^2(t)\,dt = \sum_{i,j} \lambda_{ij}^2 [\hat{\hat{Z}}_n(i,j)]^2$ where $\hat{\hat{Z}}_n(i,j) \leftrightarrow h_{ij}$ falls under the statistics covered by Theorem 3.2.

An interesting basis to consider is the set of normalized Hermite polynomials $h_j(x) = (-1)^j \cdot \frac{d^j \varphi(x)}{dx^j}/\varphi(x)$, $j \geq 3$. Here $h_3$ and $h_4$ correspond to skewness and kurtosis so that it is attractive to make $\lambda_3 = \lambda_4 = 1$ and $\lambda_j$ decrease rapidly further on.

We stress again that Propositions 4.1 and 4.2 can be applied to all these diverse tests.

4.3. *Independence.* One of the most important semiparametric hypotheses corresponds to $X = (U,V) \sim P$, $H: P = P_U \times P_V$, $U$ and $V$ are independent, $U, V \in R$. In this case the NPMLE of $P$ under $H$, known to be efficient, is $P_n = P_{nU} \times P_{nV}$, where $P_{nU}$ and $P_{nV}$ are the empirical marginals of $U$ and $V$, and is known to be efficient ([3], Chapter 5). Thus

$$\hat{\hat{Z}}_n(\gamma) = \sqrt{n}(P_n - (P_{nU} \times P_{nV}))(h)$$
$$= n^{1/2}\left\{\frac{1}{n}\sum_{i=1}^n h(U_i, V_i, \gamma, \hat{P}_n) - \frac{1}{n^2}\sum_{i=1}^n\sum_{j=1}^n h(U_i, V_j, \gamma, \hat{P}_n)\right\}.$$

Natural $h_j(u,v)$ here are of the form $h_{1\gamma}(u)h_{2\gamma}(v)$. If we take $h_\gamma(u,v) = \mathbb{1}_{Q_\gamma}(u,v) = \mathbb{1}_{Q_{1\gamma}}(u)\mathbb{1}_{Q_{2\gamma}}(v)$, where $Q_j = Q_{1\gamma} \times Q_{2\gamma}$, $Q_{j\gamma} = (-\infty, \gamma_j]$, $j = 1, 2$, we arrive at the familiar

$$\hat{\hat{Z}}_n(\gamma) = \sqrt{n}(F_n(\gamma_1, \gamma_2) - F_{nU}(\gamma_1)F_{nU}(\gamma_2)),$$

where $F_n, F_{nU}, F_{nV}$ are the appropriate empirical d.f.'s.

Application of Theorem 3.1 here is appropriate and easy. (M0) simply says $\gamma \to \mathbb{1}((x,y) \in Q_\gamma) - F_U(\gamma_1)\mathbb{1}(y \leq \gamma_2) - F_V(\gamma_2)\mathbb{1}(x \leq \gamma_1)$ is a Donsker class, essentially a statement about the bivariate empirical process. Since $h_\gamma$ does not depend on $\alpha$, (M3) and (M4) are immediate. Finally, (M5) is well known for this process—see, for instance, [3].

If we take $\mu_\alpha(d\gamma) = dF_U(\gamma_1)\,dF_V(\gamma_2)$ with $\alpha = (F_U, F_V)$ we obtain the Kiefer–Wolfowitz statistic

$$T = n\int\int (F_n(\gamma_1, \gamma_2) - F_{nU}(\gamma_1)F_{nV}(\gamma_2))^2\,dF_{nU}(\gamma_1)\,dF_{nV}(\gamma_2).$$

If we take $T = \sup_\gamma |\hat{\hat{Z}}_n(\gamma_1, \gamma_2)|$ we obtain the Kolmogorov–Smirnov version of the Kiefer–Wolfowitz statistic.

Invariance under monotone transformation of $H$ suggests

(13) $$h(u, v, \gamma, P_U \times P_V) = \mathbb{1}_{Q_j}(F_U(u), F_V(v)),$$

where $F_U, F_V$ are the c.d.f.'s of $U, V$, and leads to $\hat{\hat{Z}}_n(\gamma)$, a linear rank test statistic,

$$\hat{\hat{Z}}_n(\gamma) = n^{-1/2}\sum_{i=1}^n h_\gamma\left(\frac{R_i}{n}, \frac{S_i}{n}\right) - \frac{1}{n^{3/2}}\sum_{i=1}^n\sum_{j=1}^n h_\gamma\left(\frac{i}{n}, \frac{j}{n}\right),$$

where $R_i$ is the rank of $U_i$ among the $U$'s and $S_i$ is the rank of $V_i$ among the $V$'s. These are the building blocks of the Kallenberg–Ledwina [18] statistics, though the ones they propose are of non-$n^{-1/2}$ consistent type. We leave it to the reader to construct tests with power against all $n^{-1/2}$ alternatives and directions that (s)he prefers.

PROPOSITION 4.3. *Suppose $h_\gamma$ is given by* (13), *$\hat{\alpha}$ is the NPMLE as specified and $\alpha_0$ has continuous marginals. Then,* (M0), (M3), (M4) *and* (M5) *hold.*

PROOF. In this case

$$\Pi(h,\alpha)(x,y) = \int h(x,v)\,dF_V(v) + \int h(u,y)\,dF_U(u)$$
$$- 2\int h(u,v)\,dF_U(u)\,dF_V(v),$$

where $F_U, F_V \leftrightarrow \alpha$. So (M5) can be written

$$\sup_\gamma \left\{\left|\int h_\gamma(x,y)\,d(\mathbb{F}_{nU}(x) - F_U(x))\,d(\mathbb{F}_{nV}(y) - F_V(y))\right|\right\} = o_P(n^{-1/2}), \tag{14}$$

where $\mathbb{F}_{nU}$ is the empirical d.f. of $U$, $F_{0U}$ corresponds to $\alpha_0$, and so on. For this condition and all subsequent ones, we can assume, w.l.o.g., that $\alpha_0$ is the uniform distribution on the unit square by making separate probability integral transforms. But (14) is just

$$\sup_{0\leq\gamma_1\leq 1}|\mathbb{F}_{nU}(\gamma_1) - \gamma_1| \sup_{0\leq\gamma_2\leq 1}|\mathbb{F}_{nv}(\gamma_2) - \gamma_2| = O_P(n^{-1}),$$

and (M5) follows. (M0) has been discussed in connection with $h_\gamma$. For (M3) write

$$(P_n - P_0)(h_\gamma(\cdot,\hat{\alpha}) - h_\gamma(\cdot,\alpha_0))$$
$$= (P_n - P_0)(\mathbf{1}((u,v) \leq (F_U^{-1}(\gamma_1), F_V^{-1}(\gamma_2))) - \mathbf{1}((u,v) \leq (\gamma_1,\gamma_2))),$$

where $(\hat{F}_U, \hat{F}_V) \leftrightarrow \hat{\alpha}$. Now, by Glivenko–Cantelli $\sup|\hat{F}_U^{-1}(\gamma_1) - \gamma_1| \to 0$ and $\sup|\hat{F}_V^{-1}(\gamma_2) - \gamma_2| \to 0$. So (M3) follows from the weak convergence of $n^{1/2}(P_n - P_0)(h_\gamma)$. (M4) is argued similarly. $\square$

Application of the type I bootstrap is straightforward when $F_U$ and $F_V$ are continuous under $\alpha_0$: In view of the invariance, we need only simulation under the uniform distribution on the unit square. Resampling from the empirical (type II) is also possible but the argument is more delicate.

This result can easily be extended to more general $h_{1\gamma}(u)h_{2\gamma}(v)$ and we can also tailor tests here. For instance, consider the tensor wavelet basis on $[0,1] \times [0,1]$, $\{h_{i,j_1,j_2}\}$ where $i$ corresponds to scale and $(j_1, j_2)$ to the location. We can again suppose that departures from independence at lower resolution are more significant and proceed as in Section 4.1 to form

$$T = \int_{I^2} \Biggl( \int_{I^2} \sum_{i,j_1,j_2} \lambda_i h_{i,j_1,j_2}$$

$$\times (\mathbb{F}_{nU}(x), \mathbb{F}_{nV}(y)) \, d(P_n - \hat{P}_n)(x,y) h_{i,j_1,j_2}(u,v) \bigg)^2 \, du \, dv$$

$$= \sum_{i,j_1,j_2} \lambda_i^2 \bigg( \int h_{i,j_1,j_2}(\mathbb{F}_{nU}(x), \mathbb{F}_{nV}(y)) \, d(P_n - \hat{P}_n)(x,y) \bigg)^2,$$

where $I^2$ is the unit square and $\hat{P}_n = P_{nU} \times P_{nV}$. The $\lambda_i$ can be chosen so as to weight the lower-resolution terms as one pleases.

4.4. *Copula models.* The standard copula model is $X = (U,V)$, $U,V \in R$ as above, where for some monotone strictly increasing transformations $a(\cdot): R \to R$, $b(\cdot): R \to R$ the vector $(a(U), b(V)) \sim P_\vartheta$, $\vartheta \in \Theta$, a regular parametric model. A natural model to consider here is the bivariate Gaussian copula, where under $\vartheta$, $X$ has standard normal marginals with correlation $\vartheta$, $-1 < \vartheta < 1$. (Assuming unknown means and variances adds nothing since making $a$ and $b$ arbitrary makes these parameters unidentifiable.) In such a model we consider two problems:

(i) $H: P \in \mathcal{P}_0 = \{P_{\vartheta,a,b} : \vartheta \in \Theta, a, b \text{ general}\}$, the copula model hypothesis.

(ii) $H: \vartheta = \vartheta_0$. $K: \vartheta \neq \vartheta_0$ within $\mathcal{P}_0$.

The first hypothesis requires use of efficient estimates of $(a,b,\vartheta)$. These are in general difficult to construct. Inefficient estimates are readily computable, but application of Theorem 3.1 requires computation of $\Pi(\cdot, \alpha)$ which can be characterized by Sturm–Liouville equations and computed numerically (see [3], pages 172–175). We do not pursue this interesting special case further.

On the other hand, by our assumptions in Section 3, tests of (ii) are naturally based on $h_\gamma(\cdot, \alpha) = h_\gamma(F_U(\cdot), F_V(\cdot))$, where the $\{h_\gamma(x,y), \gamma \in R^p\}$ are scores of the parametric model $\{P_\vartheta : \theta \in \Theta\}$ at $\vartheta = \vartheta_0$. If efficient estimates $\hat{F}_U$, $\hat{F}_V$ under $H$ are used, then $\hat{\tilde{Z}}_n(\gamma)$ is the asymptotically most powerful score test in direction $\gamma$. Otherwise, if, say, we use the empirical d.f.'s $F_{nU}$ and $F_{nV}$, we can construct $\Pi(\cdot, \gamma)$, the projection on the tangent space of $\mathcal{P}_{\vartheta_0} = \{P_{(\vartheta_0,a,b)} : a,b, \text{ arbitrary}\}$ and use $\hat{Z}_n$. Finally, if efficient estimates of $\vartheta$ under $\mathcal{P}$ are available, these can be used in the obvious way, though in general such estimates will be difficult to obtain. These hypotheses of finite codimension are the subject of Choi, Hall and Schick [8]. If we specialize to the Gaussian copula model and consider $H: \rho = 0$, the independence hypothesis, it is easy to see that the single asymptotically most powerful test is to use the normal score rank statistic

$$T \equiv \frac{1}{\sqrt{n}} \sum_{i=1}^{n} \Phi^{-1}(R_i/n + 1) \Phi^{-1}(S_i/n + 1).$$

The reason here is that $F_{nU}$, $F_{nV}$ are efficient in this case—as we have already seen. Remarkably, Klaassen and Wellner [20] show that $T$ is the asymptotically most powerful score statistic for $H:\rho = \rho_0$, any $\rho_0$, by showing effectively that

$$T = \frac{1}{\sqrt{n}} \sum_{i=1}^n h(\Phi^{-1} F_U(U_i), \Phi^{-1} F_V(V_i), \rho_0) + o_P(1),$$

where $h(\Phi^{-1} F_U(U_i), \Phi^{-1} F_V(V_i), \rho_0)$ is orthogonal to the tangent space $\overset{\bullet}{\mathcal{P}}_2(\rho_0, (a_0, b_0))$, $a_0 = \Phi^{-1} F_U$, $b_0 = \Phi^{-1} F_V$.

The development of tailored tests for independence in copula models in general should be the same as for independence in general.

## APPENDIX

**Proof of Theorem 3.2.** We need

LEMMA A.1. *If* (N1) *and* (N2) *hold, then so does* (M2),

$$\sup\{\|P_{\hat{\alpha}}(h) - P_0(h) + P_0(\Pi(h, \hat{\alpha}))\|_\infty : h \in \mathcal{H}\} = o_P(n^{-1/2}).$$

PROOF. By (N1) we may w.l.o.g. assume $\hat{\alpha} \in A_0$. For simplicity let $\mu \equiv P_0$ so that $s_0 \equiv 1$. Obvious modifications suffice if $\mu \gg P_0$. Then

$$\begin{aligned}(\text{A.1})\quad P_{\hat{\alpha}}(h) - P_0(h) &= \int h\hat{s}^2 \, d\mu - \int h \, d\mu \\ &= 2\int h(\hat{s} - 1) \, d\mu + \int h(\hat{s} - 1)^2 \, d\mu.\end{aligned}$$

Since $P_{\hat{\alpha}} \Pi(h, \hat{\alpha}) = 0$,

$$\begin{aligned}(\text{A.2})\quad P_0(\Pi(h, \hat{\alpha})) &= -\int \Pi(h, \hat{\alpha})(\hat{s}^2 - 1) \, d\mu \\ &= -2\int \Pi(h, \hat{\alpha})(\hat{s} - 1) \, d\mu + \int \Pi(h, \hat{\alpha})(\hat{s} - 1)^2 \, d\mu.\end{aligned}$$

But, since $P_{\hat{\alpha}} \ll P_0$ we have by [3], formula (4b), page 50,

$$\Pi(h, \hat{\alpha}) = \hat{s}^{-1} \Pi(h\hat{s}, \alpha_0).$$

Therefore, if $\frac{\hat{s}-1}{\hat{s}} \in L_2(\mu)$,

$$\int \Pi(h, \hat{\alpha})(\hat{s} - 1) \, d\mu = \int \Pi(h\hat{s}, \alpha_0) \frac{(\hat{s} - 1)}{\hat{s}} \, d\mu$$
$$= \int h\hat{s} \Pi\left(\frac{\hat{s} - 1}{\hat{s}}, \alpha_0\right) d\mu.$$

Substituting in (A.2) we get after some manipulation

$$P_0(\Pi(h,\hat{\alpha})) = -2\int h\Pi\left(\frac{\hat{s}-1}{\hat{s}},\alpha_0\right)\hat{s}\,d\mu + \int \Pi(h,\hat{\alpha})(\hat{s}-1)^2\,d\mu$$

$$= -2\left\{\int h\Pi(\hat{s}-1,\alpha_0)\,d\mu + \int h(\hat{s}-1)\Pi(\hat{s}-1,\alpha_0)\,d\mu\right.$$

(A.3)
$$\left. - \int \hat{s}\Pi\left(\frac{(\hat{s}-1)^2}{\hat{s}},\alpha_0\right)d\mu\right\}$$

$$+ \int \Pi(h,\hat{\alpha})(\hat{s}-1)^2\,d\mu$$

$$= -2(I+II+III)+IV.$$

We bound the last three terms in absolute value by

$$|II| \le M\int |\hat{s}-1|\Pi(\hat{s}-1,\alpha_0)\,d\mu \le M\|\hat{s}-1\|_\mu^2 \le M\left\|\frac{(\hat{s}-1)^2}{\hat{s}}\right\|_\mu,$$

$$|III| \le M\left\|\Pi\left(\frac{(\hat{s}-1)^2}{\hat{s}},\alpha_0\right)\right\|_\mu \le M\left\|\frac{(\hat{s}-1)^2}{\hat{s}}\right\|_\mu,$$

where $M = \sup_{\mathcal{H}}\{\|h\|_\infty + \|\Pi(h,\alpha_0)\|_\infty\} < \infty$ by (N2).

Again, using (N2),

$$|IV| \le M\|\hat{s}-1\|_\mu^2 \le M\left\|\frac{(\hat{s}-1)^2}{\hat{s}}\right\|_\mu.$$

Combining (A.1)–(A.3) we obtain

$$P_{\hat{\alpha}}(h) - P_0(h) + P_0(\Pi(h,\hat{\alpha})) = 2\int h((\hat{s}-1) - \Pi(\hat{s}-1,\alpha_0))\,d\mu$$

$$+ O_P\left(\left\|\frac{(\hat{s}-1)^2}{\hat{s}}\right\|_\mu\right)$$

$$= O_P\left(\|\hat{s}-1-\Pi(\hat{s}-1,\alpha_0)\|_\mu + \left\|\frac{(\hat{s}-1)^2}{\hat{s}}\right\|_\mu\right)$$

$$= o_P(n^{-1/2})$$

by (N1) and the lemma follows. $\square$

PROOF OF THEOREM 3.2.  Write $\hat{h}$ for $h_\gamma(\cdot,\hat{\alpha})$ and $h$ for $h_\gamma(\cdot,\alpha_0)$. Then,

$$\hat{Z}_n(\gamma) = n^{1/2}(P_n\hat{h} - P_{\hat{\alpha}}\hat{h} - P_n\Pi(\hat{h},\hat{\alpha}))$$

$$= n^{1/2}\{(P_n - P_0)(\hat{h}) - (P_{\hat{\alpha}} - P_0)(\hat{h}) - P_n\Pi(\hat{h},\hat{\alpha})\}$$

$$= n^{1/2}(P_n - P_0)(\hat{h} - \Pi(\hat{h},\hat{\alpha})) + o_P(1)$$

by (M2). But

$$n^{1/2}(P_n - P_0)(\hat{h} - \Pi(\hat{h}, \hat{\alpha})) = n^{1/2}(P_n - P_0)(h - \Pi(h, \alpha_0)) + o_P(1)$$
$$= Z_n(\gamma, \alpha_0) + o_P(1)$$

by (M3) and (M1). (The equivalences hold uniformly in $\gamma$ by assumption.)
□

P. J. BICKEL
DEPARTMENT OF STATISTICS
UNIVERSITY OF CALIFORNIA
BERKELEY, CALIFORNIA 94720-3860
USA
E-MAIL: bickel@stat.berkeley.edu

Y. RITOV
DEPARTMENT OF STATISTICS
HEBREW UNIVERSITY OF JERUSALEM
JERUSALEM 91905
ISRAEL
E-MAIL: yaacov@mscc.huji.ac.il

T. M. STOKER
SLOAN SCHOOL OF MANAGEMENT, E52-444
MASSACHUSETTS INSTITUTE OF TECHNOLOGY
50 MEMORIAL DRIVE
CAMBRIDGE, MASSACHUSETTS 02142
USA
E-MAIL: tstoker@mit.edu